\documentclass[leqno,12pt]{article}
\usepackage{ams}
\newtheorem{definition}{Definition}[section]
\newtheorem{lemma}[definition]{Lemma}
\newtheorem{prop}[definition]{Proposition}
\newtheorem{thm}[definition]{Theorem}

\newtheorem{rmk}[definition]{Remark}

\def\LL{{\cal L}}

\def\newline{\hfil\break}

\def\dot{\hskip -.2cm {\bf .}\hskip .2cm}

\def\mapright#1{\mathop{\vbox{\ialign{
				##\crcr
    ${\scriptstyle\hfil\;\;#1\;\;\hfil}$\crcr
 \noalign{\kern-1pt\nointerlineskip}
    \rightarrowfill\crcr}}\;}}

\def\mapleft#1{\mathop{\vbox{\ialign{
				##\crcr
    ${\scriptstyle\hfil\;\;#1\;\;\hfil}$\crcr
    \noalign{\kern-1pt\nointerlineskip}
    \leftarrowfill\crcr}}\;}}

\def\proof{
  \noindent
  {\bf Proof:}
}

\def\endproof{
{\unskip\nobreak\hfil\penalty50\hskip2em\hbox{}\nobreak\hfill
          $\square$\bigbreak}
}

\usepackage{amssymb}
\begin{document}

\newcommand{\gon}{\mathop{\rm gon}\nolimits}
\newcommand{\Cliff}{\mathop{\rm Cliff}\nolimits}

\title{Green-Lazarsfeld's Gonality Conjecture for a\\
Generic Curve of Odd Genus}
\author{Marian Aprodu}
\date{}
\maketitle

\begin{abstract}
\noindent
We prove Green-Lazarsfeld's gonality conjecture for generic curves
of odd genus. The proof uses, among other things, the main result of
\cite{Ap}, and Green's syzygy canonical conjecture for generic curves of
odd genus, \cite{Vo2}. The even-genus case was previously solved in the
joint work \cite{AV}.
\end{abstract}

\section{Introduction.}

Let $C$ be a smooth, projective, irreducible, complex  curve, and $L$
a globally generated line bundle on $C$, and denote by $K_{p,q}(C,L)$
the Koszul cohomology of $C$ with value in $L$. By
definition, it is the cohomology 
of the complex:
$$
\bigwedge^{p+1} H^0(C,L)\otimes H^0(C,L^{q-1})
\rightarrow
\bigwedge^p H^0(C,L)\otimes H^0(C,L^q)
\rightarrow
\bigwedge^{p-1} H^0(C,L)\otimes H^0(C,L^{q+1}).
$$

Many geometric phenomena regarding the image of
$C$ in ${\Bbb P}H^0(C,L)^*$, such as projective
normality, or being cut out by quadrics 
can be expressed in terms of vanishing of Koszul
cohomology. 

There are a number of preliminary results which suggested that 
other geometric properties can be related to Koszul cohomology. 
A good example of a possible interaction between
geometry and Koszul cohomology is provided by the 
Nonvanishing Theorem of Green and Lazarsfeld
(cf. \cite[Appendix]{Gr1}). 
In a particular case, it shows that curves with special 
geometry have non-trivial Koszul cohomology groups. More 
precisely, if we set $d$ the gonality of $C$, and we suppose 
that $L$ is of sufficiently large degree, then one has:
$$
K_{h^0(C,L)-d-1,1}(C,L)\not= 0.
$$ 

It was conjectured by Green and Lazarsfeld, 
\cite[Conjecture 3.7]{GL1}, that this
result were optimal, that is:

\begin{equation}
\label{equation: gonality conjecture}
K_{h^0(C,L)-d,1}(C,L)= 0,
\end{equation}
\\
for all $L$ of sufficiently large degree. Broadly speaking,
it would mean that gonality of curves could be read off 
minimal resolutions of embeddings of sufficiently large degree.

We prove here the following:

\begin{thm}\dot
\label{thm: main}
Green-Lazarsfeld's gonality conjecture holds 
for a generic curve of odd genus.
\end{thm}

The case of generic curves of odd genus is exactly the case
left away from the joint work \cite{AV}, where the gonality
conjecture was verified for a generic curve of genus $g$ and 
gonality $d$ with $\frac{g}{3} + 1\leq d < [\frac{g + 3}{2}]$.
Recall that curves of given genus $g$ and gonality $d$
are parametrized by an irreducible quasi-projective
variety. Beside, we always have $d \leq [\frac{g + 3}{2}]$,
and the maximal gonality is realized on a non-empty open
set of moduli space of curves of genus $g$. 
Therefore, Theorem \ref{thm: main} and the results of
\cite{AV} complete each other, showing eventually
that Green-Lazarsfeld's gonality conjecture holds
for a generic curve of genus $g$ and gonality
$d$ with $\frac{g}{3} + 1\leq d$.

The techniques used for proving Theorem \ref{thm: main} 
are standard, and very similar to those used in \cite{Ap}, 
\cite{AV}, \cite{Vo2}, see Section 2 for details. Green's
conjecture plays a central role in the whole proof.
An alternative proof of the generic gonality conjecture 
for the even-genus case is sketched in Remark \ref{rmk: even}. 
In Section 3, we prove a result, Proposition 
\ref{prop: two points}, which makes a clear difference between the 
case of generic curves of odd genus,
and the curves analysed in \cite{AV}, justifying
the strategy chosen for proving Theorem \ref{thm: main}.

\bigskip

\noindent
{\em Notation and conventions.}
\\\\
If $V$ a finite-dimensional complex vector space, and $S(V)$ 
denotes its symmetric algebra, for any graded
$S(V)$--module
$B=\mathop\bigoplus\limits_{q\in {\Bbb Z}}B_q$
there is a naturally defined
complex of vector spaces, called the {\em Koszul complex} of $B$
(cf. \cite[Definition 1.a.7]{Gr1}),
$$
...\longrightarrow 
\bigwedge ^{p+1}V\otimes B_{q-1}
\stackrel{d_{p+1,q-1}}{\longrightarrow}
\bigwedge ^pV\otimes B_q
\stackrel{d_{p,q}}{\longrightarrow}
\bigwedge ^{p-1}V\otimes B_{q+1}
\longrightarrow ... ,$$
where $p$, and $q$ are integer numbers.
The cohomology of this complex is denoted by 
$$ K_{p,q}(B,V)=\mbox{Ker }d_{p,q}/ 
\mbox{Im }d_{p+1,q-1}. $$

In the algebro-geometric context, if
$X$ is a complex projective variety, 
$L\in \mbox{Pic}(X)$ is a line
bundle, $\cal F$ is a coherent sheaf,
and $V=H^0(X,L)$, one usually computes 
Koszul cohomology for the $S(V)$-module
$B=\bigoplus\limits_{q\in{\Bbb Z}}
H^0(X,L^q\otimes{\cal F})$.
The standard notation is $K_{p,q}(X,{\cal F},L) =K_{p,q}(B,V)$. 
If ${\cal F}\cong {\cal O}_X$, we
drop it, and write simply $K_{p,q}(X,L)$.

\bigskip

Throughout this paper, we shall refer very often to \cite{Gr1}
and \cite{Gr2} for the basic facts of the Koszul cohomology theory.

\section{Proof of Theorem \ref{thm: main}.}

A generic curve of odd genus $2k+1$ is of gonality
$k+2$. Since elliptic curves are fairly well
understood, we can suppose $k\geq 1$.
Recall that vanishing of Koszul cohomology is an open 
property (see, for example, \cite{BG}), and the moduli 
space of curves of given genus is irreducible, so it 
would suffice to exhibit one curve $C$ of odd genus $2k+1$, 
and one nonspecial line bundle $L_C$ on $C$ such that 
$$
K_{h^0(C,L_C)-k-2,1}(C,L_C)=0.
$$

This principle has already been used in \cite{Ap} and \cite{AV},
and is based on the main result of \cite{Ap}
(see also \cite[Theorem 2.1]{AV}):

\begin{thm}\dot
\label{thm: Ap}
If $L$ is a nonspecial line bundle on a curve 
$C$, which satisfies\linebreak
$K_{n,1}(C,L)=0$, for a positive integer
$n$, then, for any effective divisor $E$ of degree $e$,
we have $K_{n+e,1}(C,L+E)=0$.
\end{thm}

We construct $C$ and $L_C$ in the following way.
Let $S$ be a $K3$ surface whose Picard group is 
freely generated by a very ample line bundle $\LL$, 
and by one smooth rational curve $\Delta$, 
such that $\LL ^2=4k$, and $\LL.\Delta= 3$.
Such a surface does exist, as one can see by 
analysing the period map (see, for example,
\cite[Lemma 1.2]{Og}). Denote $L=\LL +\Delta$.
Then one can show:

\begin{lemma}\dot
\label{lemma: 1}
All smooth curves in the linear systems
$|\LL|$ and $|L|$ have maximal Clifford index,
and thus maximal gonality, too.
\end{lemma}

\proof
The proof runs exactly like in
\cite[Proposition 1]{Vo2},
using \cite{GL2}.

\endproof

Let $C\in |\LL|$ be a smooth curve, and set
$L_C=L_{|C}$. Then $C$ is of genus $g=2k+1$, 
and $L_C$ is of degree $2g+1=4k+3$, thus nonspecial.
Beside, the curve $C$ is generic, and
$h^0(C,L_C)=2k+3$, so the vanishing 
(\ref{equation: gonality conjecture})
predicted by the gonality conjecture for
the pair $(C,L_C)$ becomes:

\begin{equation}
\label{equation: desired vanishing}
K_{k+1,1}(C,L_C)=0.
\end{equation}

This vanishing will be a consequence of
another two Lemmas which are proven below. 

\bigskip

\begin{lemma}\dot
\label{lemma: 2}
$K_{k+1,1}(S,L)=0$.
\end{lemma}

\proof
Remark that any smooth curve $D\in|L|$ is of genus $2k+3$.
From \cite{Vo2}, \cite{HR}, and Lemma \ref{lemma: 1}, it follows that
Green's conjecture is valid for $D$. Therefore
$K_{k+1,1}(D,K_D)=0$. From Green's hyperplane section
theorem \cite[Theorem 3.b.7]{Gr1}, we conclude $K_{k+1,1}(S,L)=0$.

\endproof

\begin{lemma}\dot
\label{lemma: 3}
For any integer $p$, we have a natural 
isomorphism\linebreak
$K_{p,1}(S,L)\cong K_{p,1}(C,L_C)$.
\end{lemma}

\proof
We set $D=C+\Delta$; it is a connected reduced
divisor on $S$. All the hypotheses of  Green's hyperplane section
theorem \cite[Theorem 3.b.7]{Gr1} are fulfilled, so
one has isomorphisms
$$
K_{p,1}(S,L)\cong K_{p,1}(D,L_D),
$$
for all integers $p$, where we denoted $L_D=L_{|D}$.

We use next the natural short exact sequence
of coherent sheaves on $S$ (see also \cite{AN}):
$$
0\rightarrow {\cal O}_\Delta (-C)
\rightarrow {\cal O}_D \rightarrow {\cal O}_C
\rightarrow 0,
$$
which yields furthermore to an exact sequence, 
for any integer $q$:
$$
0\rightarrow H^0(\Delta,{\cal O}_\Delta (L^{q-1}+\Delta))
\rightarrow H^0(D,L^q_D) \rightarrow H^0(C,L^q_C)
$$

We analyse the maps $H^0(D,L^q_D) \rightarrow H^0(C,L^q_C)$
for different $q$. If $q=0$ the corresponding map is obviously
an isomorphism. For $q=1$, we also have an isomorphism
$H^0(D,L_D) \cong H^0(C,L_C)$,
since ${\cal O}_\Delta (\Delta)\cong K_\Delta 
\cong {\cal O}_{{\Bbb P}^1}(-2)$.
Moreover, since $(L+\Delta).\Delta <0$, for $q=2$, we
obtain an inclusion $H^0(D,L_D^2)\subset H^0(C,L_C^2)$.
These facts reflect into an isomorphism
between Koszul cohomology groups
$K_{p,1}(D,L_D)\cong K_{p,1}(C,L_C)$ for any
integer $p$ (apply, for example,
\cite[Remark 1.1]{Ap}). Therefore, for any $p$,
$K_{p,1}(S,L)$ is isomorphic to
$K_{p,1}(C,L_C)$.

\endproof

Summing up, $K_{k+1,1}(C,L_C)$ is isomorphic to
$K_{k+1,1}(S,L)$, and the latter vanishes,
which proves (\ref{equation: desired vanishing}),
and Theorem \ref{thm: main}, too.

\endproof

\begin{rmk}\dot
\label{rmk: even}
{\rm
A similar idea can be used to give an
alternative proof of the gonality conjecture
for a generic curve of even genus. Let us
sketch the proof in a few words.

\bigskip

Let $S$ be a $K3$ surface whose Picard group is
freely generated by a very ample line bundle $\LL$,
and by one smooth rational curve $\Delta$,
such that $\LL ^2=4k-2$, and $\LL.\Delta= 2$, and
let furthermore $C\in |\LL|$ be
a smooth curve. Denote $L=\LL +\Delta$, and
$L_C=L_{|C}$.
Then $C$ is of genus $g=2k$ and gonality $k+1$, 
and $L_C$ is of degree $2g$.

Since $h^0(C,L_C)=2k+1$, the vanishing
(\ref{equation: gonality conjecture})
predicted by the gonality conjecture for
the pair $(C,L_C)$ becomes
$K_{k,1}(C,L_C)=0$.

A smooth curve in the linear system $|L|$ has
genus $2k+1$ and gonality $k+2$, and thus it 
satisfies Green's conjecture, \cite{Vo2}.
It implies $K_{k,1}(S,L)=0$, after having applied
Green's hyperplane section theorem.

The proof of Lemma \ref{lemma: 3} can be easily
addapted to obtain an isomorphism $K_{p,1}(S,L)\cong
K_{p,1}(C,L_C)$, for any $p$. In particular, for
$p=k$, we obtain $K_{k,1}(C,L_C)=0$.  
}
\end{rmk}

\begin{rmk}\dot
{\rm 
It follows from Theorem \ref{thm: main}
and Theorem \ref{thm: Ap} that the
vanishing (\ref{equation: gonality conjecture})
is valid for any line bundle of degree at least
$3g+1$ over a generic curve of odd genus $g$.
In the even-genus case, the lower bound
was $3g$, see \cite{AV}. We could predict thus that, up 
to one unit, the lower bound $3g$, indicated 
by Green and Lazarsfeld as one of the possible 
replacements for the more ambiguous 
"{\em of sufficiently large degree}", 
is correct.
}
\end{rmk}

\section{Why adding three points to the canonical bundle?}

The strategy used in \cite{AV} to verify the gonality 
conjecture for a curve $C$ of genus $g$, and gonality 
$d$ with $\frac{g}{3} + 1\leq d < [\frac{g + 3}{2}]$
was to look at bundles of type $K_C+x+y$, where $x$, 
and $y$ were suitably chosen points of $C$ (see also 
Remark \ref{rmk: even}).
For Theorem \ref{thm: main} instead, we needed to add three
points to the canonical bundle in order to have the desired 
vanishing (\ref{equation: gonality conjecture}).
This choice is justified by the following remark, which
shows that on a generic curve $C$ of odd genus, bundles of type
$K_C+x+y$ never verify (\ref{equation: gonality conjecture}).

\begin{prop}\dot
\label{prop: two points}
Let $C$ be a curve of odd genus
$2k+1\geq 3$, and maximal gonality $k+2$. 
Then, for any two points $x$, and $y$
of $C$, the dimension of
$
K_{k,1}(C,K_C+x+y)
$
equals the binomial coefficient
$\left(  
\begin{array}{c}
2k+1 \\
k+2
\end{array}
\right)$.
\end{prop}

For the proof of Proposition \ref{prop: two points}, 
we need the following elementary Lemma.

\begin{lemma}\dot
\label{lemma: multiplication}
Let $X$ be an irreducible projective manifold,
and $D\not=0$ be an effective divisor.
Then, for any $L\in\mbox{\rm Pic}(X)$, and any integer $p\geq 1$,
we have an exact sequence:
$$
0\rightarrow \bigwedge ^{p+1}H^0(X,L-D)\rightarrow K_{p,1}(X,-D,L)
\rightarrow K_{p,1}(X,L).
$$ 
\end{lemma}

\proof
Let $V=H^0(X,L)$, and consider the graded $S(V)$--modules
$$
A=\mathop\bigoplus\limits_{q\in{\Bbb Z}} H^0(X,L^q-D),\;
B=\mathop\bigoplus\limits_{q\in{\Bbb Z}} H^0(X,L^q),
$$
and $C=B/A$, where the inclusion of $A$ in $B$
is given by the multiplication with the non-zero section
of ${\cal O}_X(D)$ vanishing along $D$.
Obviously, $A_0=0$, and $C_0\cong {\Bbb C}$.
The long cohomology sequence for syzygies 
yields to an exact sequence:

$$
0\rightarrow
\mbox{Ker}\left( \bigwedge ^{p+1}V\rightarrow \bigwedge ^pV
\otimes C_1\right)\rightarrow K_{p,1}(X,-D,L)\rightarrow
K_{p,1}(X,L)\rightarrow ...
$$

We aim to prove that 

$$
\mbox{Ker}\left( \bigwedge ^{p+1}
V\rightarrow\bigwedge ^pV \otimes C_1\right)\cong
\bigwedge ^{p+1}H^0(X,L-D).
$$

For this, choose a basis $\{ w_1,...,w_N\}\subset V$,
such that $\{ w_1,...,w_s\}$ is a basis of
$H^0(X,L-D)$, and pick an element
$\alpha =\sum\limits_{1\leq i_1<...<i_{p+1}\leq N}
\alpha_{i_1...i_{p+1}}
w_{i_1}\wedge...\wedge w_{i_{p+1}}\in\bigwedge ^{p+1}V$.
It belongs to $\mbox{Ker}\left( \bigwedge ^{p+1}
V\rightarrow\bigwedge ^pV \otimes C_1\right)$ if
and only if the following relations are satisfied, for
any $1\leq k_1<...<k_p\leq N$:

$$
\sum\limits_{k\not\in\{k_1,...,k_p\}}
(-1)^{\#\{ k_i<k\}}\alpha_{k_1...k...k_p}w_k\in
H^0(X,L-D).
$$

In particular, for any $1\leq k_1<...<k_p\leq N$,
and $k>s$, $\alpha_{k_1...k...k_p}=0$,
in other words all $\alpha_{i_1...i_{p+1}}$
with $i_{p+1}>s$ vanish, so $\alpha $
belongs to $\bigwedge ^{p+1}H^0(X,L-D)$.

\endproof

\bigskip

\noindent
{\bf Proof of Proposition \ref{prop: two points}:}
We compute first the dimension of $K_{k-1,2}(C,K_C+x+y)$.
By Green's duality Theorem it equals the dimension
of $K_{k+1,1}(C,-x-y,K_C+x+y)$. We know 
$K_{k,1}(C,K_C)=0$, as $C$ satisfies Green's conjecture.
From \cite[Theorem 3]{Ap}, it follows $K_{k+1,1}(C,K_C+x+y)=0$. 
Then from Lemma \ref{lemma: multiplication} 
we obtain an isomorphism
$$
\bigwedge ^{k+2}H^0(C,K_C)\cong K_{k+1,1}(C,-x-y,K_C+x+y).
$$
Therefore, the dimension of $K_{k+1,1}(C,-x-y,K_C+x+y)$,
and thus of $K_{k-1,2}(C,K_C+x+y)$,
equals the binomial coefficient 
$\left(
\begin{array}{c}
2k+1 \\
k+2
\end{array}
\right)$. 
What is left from the proof is a combinatorial
computation. Analysing the Koszul complex which
computes $K_{k,1}(C,K_C+x+y)$, one can prove:

\begin{lemma}\dot
\label{lemma: combinatorics}
The Euler characteristic of the complex
$$
0\rightarrow \bigwedge^{k+1}H^0(K_C+x+y)
\rightarrow \bigwedge^kH^0(K_C+x+y)\otimes H^0(K_C+x+y)
\rightarrow ...
$$
equals zero.
\end{lemma}

\proof
Standard combinatorics. Use 
$h^0(C,K_C^q+qx+qy)=4qk+2q-2k$, for all $q\geq 1$.

\endproof

\bigskip

Knowing that $K_{k-j+1,j}(C,K_C+x+y)$ vanishes for
all $j\not=1,2$, we conclude
that $K_{k,1}(C,K_C+x+y)$ and 
$K_{k-1,2}(C,K_C+x+y)$ have the same dimension.

\endproof

\noindent
{\sc Author's Adresses}
\\\\
Romanian Academy, Institute of Mathematics
"Simion Stoilow", P.O.Box 1-764, RO-70700, 
Bucharest, Romania (e-mail:
Marian.Aprodu\char64 imar.ro)
\\\\
Universit\'e de Grenoble 1,
Laboratoire de Math\'ematiques,
Institut Fourier BP 74,
F-38402 Saint Martin d'H\`eres Cedex,
France (e-mail: Marian.Aprodu\char64 ujf-grenoble.fr)


\begin{thebibliography}{99}
\bibitem{Ap} M. Aprodu, On the vanishing of higher syzygies of
curves. Math. Z. {\bf 241} (2002) 1-15.
\bibitem{AN} M. Aprodu, J. Nagel, A Lefschetz type result
for Koszul cohomology. to appear
in Manuscr. Math.
\bibitem{AV} M. Aprodu, C. Voisin, 
Green-Lazarsfeld's conjecture for generic curves of large
gonality. C.R.A.S. {\bf 36} (2003) 335-339.
\bibitem{BG} M. Boraty\'nsky and S. Greco,
Hilbert functions and Betti numbers in a flat family,
Ann. Mat. Pura Appl. (4) {\bf 142} (1985)
277-292.
\bibitem{Gr1} M. Green, Koszul cohomology and the geometry of
projective varieties. J. Diff. Geom. {\bf 19} (1984) 125-171
(with an Appendix by M. Green and R. Lazarsfeld).
\bibitem{Gr2} M. Green, Koszul cohomology and the geometry of
projective varieties. II. J. Diff. Geom. {\bf 20} (1984) 279-289. 
\bibitem{GL1} M. Green and R. Lazarsfeld, On the projective normality
of complete linear series on an algebraic curve, Invent. Math. {\bf 83}
(1986) 73-90.
\bibitem{GL2} M. Green and R. Lazarsfeld, Special divisors on curves on a
$K3$ surface. Invent. Math. {\bf 89} (1987) 357-370. 
\bibitem{HR} A. Hirschowitz, Ramanan, New evidence for Green's conjecture
on syzygies of canonical curves. 
Ann. Sci. \'Ecole Norm. Sup. (4) {\bf 31} (1998) 145-152.
\bibitem{Og} K. Oguiso, Two remarks on Calabi-Yau Moishezon
threefolds. J. reine angew. Math. {\bf 452} (1994) 153-162.
\bibitem{Vo1} C. Voisin, Green's generic syzygy conjecture for curves of
even genus lying on a $K3$ surface. J. Eur. Math. Soc. {\bf 4}
(2002) 363-404.  
\bibitem{Vo2} C. Voisin, Green's canonical syzygy
conjecture for curves of odd genus, arXiv:math.AG/0301359.
\end{thebibliography}
\end{document}